\documentclass[preprint,12pt]{amsart}

\usepackage{amsfonts,amsmath,amssymb,amsxtra,amsthm}
\usepackage{mathrsfs}
\usepackage{times}
\usepackage{anysize}
\usepackage{graphicx}

\newcommand{\CM}{\mathcal{M}}
\newcommand{\XT}{\widetilde{X}}
\newcommand{\yy}{\mathbf{y}}
\newcommand{\kk}{\underline{k}}

\DeclareMathOperator{\Hh}{H}
\DeclareMathOperator{\Tr}{Tr}
\DeclareMathOperator{\Aut}{Aut}
\DeclareMathOperator{\GCD}{GCD}

\newtheorem{theorem}{Theorem}[section]
\newtheorem{lemma}[theorem]{Lemma}

\theoremstyle{definition}
\newtheorem{definition}[theorem]{Definition}
\newtheorem{example}[theorem]{Example}

\newtheorem{corollary}[theorem]{Corollary}
\newtheorem{proposition}[theorem]{Proposition}

\theoremstyle{remark}
\newtheorem{remark}[theorem]{Remark}

\title{The equivariant Euler characteristic of moduli spaces of  curves.}
\author{Eugene Gorsky}
\address{Mathematics Department, Stony Brook University, 
Stony Brook NY, 11794-3651, USA}
\email{ egorsky@math.sunysb.edu}

\begin{document}

\maketitle

\begin{abstract}
We derive a formula for the $S_n$--equivariant Euler characteristic of the moduli space $\CM_{g,n}$ of genus $g$ curves with $n$ marked points.  
\end{abstract}

\section{Introduction}

Consider the moduli space $\CM_{g,n}$ of algebraic curves of genus $g$ with $n$ marked points. The symmetric group $S_n$ acts naturally on this space. Let $V_{\lambda}$ denote the irreducible representation of $S_n$ corresponding to a Young diagram $\lambda$,  then one can decompose the cohomology of $\CM_{g,n}$ into isotypic components:
 $$\Hh^{i}(\CM_{g,n})=\bigoplus_{\lambda}a_{i,\lambda}V_{\lambda}.$$ The $S_n$-equivariant Euler characteristic of $\CM_{g,n}$ is defined by the formula
$$\chi^{S_n}(\CM_{g,n})=\sum_{i,\lambda}(-1)^{i}a_{i,\lambda}s_{\lambda},$$
where $s_{\lambda}$ denotes the Schur polynomial labeled by the diagram $\lambda$. We calculate these equivariant Euler characteristics for all $g\ge 2$ and $n$.
 
\begin{theorem}
\label{main}
The generating function for the $S_n$-equivariant Euler characteristics of $\CM_{g,n}$ has the form
$$
\sum_{n=0}^{\infty}t^{n}\chi^{S_n}(\CM_{g,n})=\sum_{\underline{k}}c_{k_1,\ldots,k_{r}}
\prod_{j=1}^{r}(1+p_jt^j)^{k_j},
$$
where $p_j$ are power sums  and the coefficients
$c_{k_1,\ldots,k_{r}}$ are defined by the equation \eqref{defC}.
\end{theorem} 

Consider the moduli space $\CM_{g}(k_1\ldots,k_r)$ of pairs $(C,\tau)$ where $C$ is a genus $g$ curve and $\tau$ is an automorphism of $C$ such that for all $i$ the Euler characteristic of the set of points in $C$ having the orbit of length $i$  under the action of $\tau$ equals $ik_i$.
The coefficient $c_{k_1,\ldots,k_{r}}$ can be also defined as the orbifold Euler characteristic of  $\CM_{g}(k_1\ldots,k_r)$. 

This moduli space can be defined for any tuple of integers $(k_1,\ldots,k_{r})$ of arbitrary size $r$, but we prove that (for a fixed genus $g$) it is non-empty only for a finite number of such tuples. In particular, $r$ cannot exceed $4g+2$.

\begin{corollary}
The generating function $\sum_{n=0}^{\infty}t^n\chi^{S_{n}}(\CM_{g,n})$ is a rational function in $t$. Furthermore,  for any $n$,
$$
\chi^{S_{n}}(\CM_{g,n})\in \mathbb{Z}[p_1,\ldots,p_{4g+2}].
$$
\end{corollary}

Theorem \ref{main} can be compared with the computations of \cite{bgp}, \cite{faber}, \cite{getzler2} and \cite{my} in genus 2 and with the computations of \cite{bergstrom}, \cite{bertom}, \cite{gelo}, \cite{tommasi} and \cite{t2} in genus 3. A similar generating function for the moduli spaces of hyperelliptic curves was previously obtained in \cite{myhyp}.
The non-equivariant Euler characteristics of moduli spaces of curves were computed by Bini and Harer in \cite{binihar}.

The paper is organized as follows. In Section \ref{sec:euchar}
we consider a complex quasi-projective variety $X$ with an action of a finite group $G$.  Theorem \ref{th:conf} provides a formula for the $S_n$-equivariant Euler characteristic of quotients $F(X,n)/G$, where $F(X,n)$ is a configuration space of $n$ labeled distinct points on $X$. This theorem  was previously proved in \cite{my} using the results of Getzler \cite{getzler,getzler1} concerning Adams operations over the equivariant motivic rings (see also \cite{myadams}). The alternative proof presented here uses only the basic properties of  Euler characteristic and seems to be more  geometric. It also makes the proof of the main result self-contained.

In Section \ref{sec:moduli} we apply this theorem to the universal family over $\CM_{g}$, the moduli space of genus $g$ curves. This allows us to prove in Theorem \ref{th:EqToOrb} that the coefficient 
$c_{k_1,\ldots,k_{r}}$ is equal to the orbifold Euler characteristic of  $\CM_{g}(k_1\ldots,k_r)$. These Euler characteristics are then computed in Theorem \ref{fin} using the results of Harer and Zagier. 

\section*{Acknowledgements}

The author is grateful to J. Bergstr\"om, S. Gusein-Zade, M. Kazaryan and S. Lando for useful discussions. This work was partially supported by the grants RFBR-007-00593, RFBR-08-01-00110-a, NSh-709.2008.1 and the M\"obius Contest fellowship for young scientists.
 
\section{Equivariant Euler characteristics}
\label{sec:euchar}

Let $X$ be a complex quasi-projective variety with an action of a finite group $G$. Let us denote by $F(X,n)$ the configuration space of ordered 
$n$-tuples of distinct points on $X$. 
For each $n$, the action of the group $G$ on $X$ can be naturally extended to the action of $G$ on $F(X,n)$, commuting with the natural action of $S_n$. 

In the computations below we will use the additivity and multiplicativity of the Euler characteristic, as well as the Fubini formula for the integration with respect to the Euler characteristic (\cite{kho,viro}, see also \cite{mcph}).

\begin{lemma}
\label{lem:power}
The following equation holds:
$\sum_{n=0}^{\infty}\frac{t^n}{n!}\chi(F(X,n))=(1+t)^{\chi(X)}.$
\end{lemma}

\begin{proof}
The map $\pi_n:F(X,n)\rightarrow F(X,n-1)$, which forgets the last point in the $n$-tuple, has fibers isomorphic to $X$ without $n-1$ points. Therefore $\chi(F(X,n))=(\chi(X)-n+1)\cdot\chi(F(X,n-1))$, and 
$\chi(F(X,n))=\chi(X)\cdot(\chi(X)-1)\cdot\ldots\cdot(\chi(X)-n+1).$
\end{proof}

Let $p_k$ denote the $k$th power sum and let $V_{\lambda}$ denote the irreducible representation of $S_n$ labelled by the Young diagram $\lambda$. We define the $S_n$-equivariant Euler characteristic of $F(X,n)/G$ by the equation
$$\chi^{S_n}(F(X,n)/G)=\sum_{i,\lambda}(-1)^{i}a_{i,\lambda}s_{\lambda},$$
where $\Hh^{i}(F(X,n)/G)=\bigoplus_{\lambda}a_{i,\lambda}V_{\lambda}$ and $s_{\lambda}$ is the Schur polynomial.

\begin{lemma}
\label{lem:frob}
The following equation holds:
$$
\chi^{S_n}(F(X,n)/G)=
\frac{1}{n!}\sum_{\sigma\in S_{n}} p_1^{k_1(\sigma)}\cdot\ldots\cdot p_{n}^{k_n(\sigma)}\cdot \chi\left([F(X,n)/G]^{\sigma}\right),
$$
where $k_i(\sigma)$ is the number of cycles of length $i$ in a permutation $\sigma$.
\end{lemma}
\begin{proof}
It is well known that for every $i$ 
$$
\sum_{\lambda} a_{i,\lambda}s_{\lambda}=\frac{1}{n!} \sum_{\sigma\in S_{n}}p_1^{k_1(\sigma)}\cdot\ldots\cdot p_{n}^{k_n(\sigma)}\cdot \Tr(\sigma)|_{\Hh^{i}(F(X,n)/G)},
$$
hence
$$
\chi^{S_n}(F(X,n)/G)=\frac{1}{n!}\sum_{i}(-1)^{i}\sum_{\sigma\in S_{n}} p_1^{k_1(\sigma)}\cdot\ldots\cdot p_{n}^{k_n(\sigma)}\cdot \Tr(\sigma)|_{\Hh^{i}(F(X,n)/G)}
$$
Now the statement follows from the Lefschetz fixed point theorem.
\end{proof}

\begin{lemma}
\label{lem:fixpt}
Let $\sigma\in S_n$. Then
$$\chi\left([F(X,n)/G]^{\sigma}\right)=\frac{1}{|G|}\sum_{g\in G}\chi\left(F(X,n)^{g^{-1}\sigma}\right).$$
\end{lemma}

\begin{proof}
For a point $\yy\in F(X,n)$ whose projection on $F(X,n)/G$ is $\sigma$-invariant there exists an element $g\in G$ such that $\sigma{\yy}=g{\yy}$. Consider the set of pairs 
$$S=\{(g,{\yy})|g\in G,{\yy}\in F(X,n)|\sigma{\yy}=g{\yy}\}$$
and its two-step projection $S\rightarrow F(X,n)\rightarrow F(X,n)/G$. The fiber of the first projection over a point ${\yy}$ is isomorphic to $G$-stabiliser of $\yy$ or empty, the fiber of the second projection containing $\yy$ is exactly the orbit of $\yy$. Therefore the cardinality of every fiber of the composition is equal to $|G|$.  
\end{proof}

\begin{definition}
For any $g\in G$ we denote by $X_k(g)$ the subset of $X$ consisting of points with $g$-orbits of length $k$. For example, $X_1(g)$ is a set of $g$-fixed points. Let $\XT_{k}(g)=X_k(g)/(g)$, where $(g)$ is
a cyclic subgroup in $G$ generated by $g$.
\end{definition}

The following theorem was deduced in \cite{my} from the results of Getzler \cite{getzler,getzler1}, here we would like to present a more geometric and straightforward proof.

\begin{theorem}
\label{th:conf}
The generating function for the $S_n$-equivariant Euler characteristics of the quotients $F(X,n)/G$ is given by the following equation:

\begin{equation}
\label{eqF}
\sum_{n=0}^{\infty}t^{n}\chi^{S_n}(F(X,n)/G)=\frac{1}{|G|}\sum_{g\in G}\prod_{k=1}^{\infty}(1+p_{k}t^k)^{\frac{\chi(X_k(g))}{k}}.
\end{equation}
\end{theorem}

\begin{proof}
Since all points in $X_k(g)$ have $g$-orbit of length $k$, we have $\chi(\XT_k(g))=\chi(X_k(g))/k.$ 
From Lemma \ref{lem:power} one gets:
$$
(1+p_{j}t^j)^{\chi\left(\XT_j(g)\right)}=\sum_{k_j=0}^{\infty}\frac{p_j^{k_j}t^{jk_j}}{(k_j)!}\chi\left(F\left(\XT_j(g),k_j\right)\right),
$$
Therefore the coefficient at $t^n$ in the right hand side of (\ref{eqF}) equals to:
$$
\frac{1}{|G|}\sum_{g\in G}\sum_{\sum jk_j=n}\prod_{j}\frac{p_j^{k_j}}{k_j!}\chi\left(F\left(\XT_j(g),k_j\right)\right).
$$
On the other hand,  by Lemma \ref{lem:frob} and Lemma \ref{lem:fixpt}, the left hand side of (\ref{eqF}) can be rewritten as following:
$$
\frac{1}{|G|}\sum_{g\in G}\frac{1}{n!}\sum_{\sigma\in S_{n}} p_1^{k_1(\sigma)}\cdot\ldots\cdot p_{n}^{k_n(\sigma)}\cdot \chi([F(X,n)]^{g^{-1}\sigma}).
$$
If for a tuple $\yy\in F(X,n)$ we have $\sigma(\yy)=g(\yy)$, the action of $(g)$ at this tuple has $k_j(\sigma)$ cycles of length $j$. Every cycle of length $j$ corresponds to a point in $\XT_j(g)$, hence for every $g$ we can define a map 
$$\alpha_{g}: \sqcup_{\sigma\in S_n} [F(X,n)]^{g^{-1}\sigma}\to \prod_{j} F(\XT_j(g),k_j)/S_{k_j}.$$
Given a $g$-invariant $n$-tuple of distinct points in $X$, there are $n!$ ways to label them and make an ordered tuple $\yy$. Every such ordering defines a unique permutation $\sigma$ such that
$\sigma (\yy)=g(\yy),$ therefore all fibers of $\alpha_{g}$
have cardinality $n!$ and
$$\frac{1}{n!}\sum_{\sigma\in S_{n}}\chi([F(X,n)]^{g^{-1}\sigma})=\prod_{j} \chi\left(F(\XT_j(g),k_j)/S_{k_j}\right)=\prod_{j} \frac{\chi\left(F(\XT_j(g),k_j)\right)}{k_j!}.$$
\end{proof}

\section{Moduli spaces of curves}
\label{sec:moduli}

Let us apply Theorem \ref{th:conf} to the study of moduli spaces of curves. 
Let $\CM_{g}$ denote the moduli space of genus $g$ algebraic curves and let $\CM_{g,n}$ denote the moduli space of genus $g$ algebraic curves with $n$ parked points (we will always assume $g\ge 2$). Let $\CM_{g}(k_1,\ldots,k_r)$ be the moduli space of pairs $(C,\tau)$ where $C$ is a genus $g$ curve and $\tau$ is an automorphism of $C$ such that $\chi(C_i(\tau))=ik_i$ for all $i$. Since $g\ge 2$, every automorphism of $C$ has finite order, hence one can choose $r$         
 such that $k_r\neq 0$ and $k_i=0$ for $i>r$.

There is a natural forgetful map $\pi_{g,\kk}:\CM_{g}(k_1,\ldots,k_r)\to \CM_{g}$ sending $(C,\tau)$ to $C$.
For a curve $C$ we define $\Aut_{\kk}(C)=\pi_{g,\kk}^{-1}(C)\subset \Aut(C).$
%Note that $\sum_{i=1}^{n}ik_i=\chi(C)=2-2g.$ 

\begin{proposition}
\label{prop:bounds}
Suppose that $\CM_{g}(k_1,\ldots,k_r)$ is not empty. 
Then $k_r<0, k_i=0$ for $i\nmid r$ and $k_i\ge 0$ for $i\mid r, i<r$. Moreover, we have the following bounds 
on $r$ and $k_i$:
$$
r\le 4g+2,\ |k_r|\le 2g,\ \sum_{i=1}^{r-1}k_i\le 2g+2.
$$ 
\end{proposition}

\begin{proof}
Let $\tau$ be an automorphism of a genus $g$ curve $C$ such that $\chi(C_{i}(\tau))=ik_i$ for all $i$. Note that $C_{i}(\tau)$ are finite sets for $i<r$ and 
\begin{equation}
\label{eq:chiC}
\chi(C)=2-2g=\sum_{i=1}^{r-1}ik_i-r|k_r|
\end{equation}
The quotient $C_1=C/\tau$  is a smooth curve of some genus $h$, and the Riemann-Hurwitz formula   yields its Euler characteristic:
\begin{equation}
\label{eq:chiC1}
\chi(C_1)=2-2h=\sum_{i=1}^{r-1}k_i-|k_r|.
\end{equation}
The projection of $C$ to $C_1$ is a ramified covering of order $r$ with $s=\sum_{i=1}^{r-1}k_j$ ramification points. The automorphism $\tau$ has order $r$, so $i|r$, if $k_i\neq 0$. By a theorem of Wiman (\cite{wiman}, see also \cite{harvey}), the maximal order for an automorphism of a genus $g$ curve equals $4g+2$, hence $r\le 4g+2$.

Since proper divisors of $r$ cannot exceed $r/2$,  equation \eqref{eq:chiC1} implies:
$$\sum_{i=1}^{r-1}ik_i\le  \frac{r}{2}  \sum_{i=1}^{r-1}k_i=\frac{r}{2}(2-2h+|k_r|),$$
hence by \eqref{eq:chiC}:
\begin{equation}
\label{eq:bound}
2g-2=r|k_r|-\sum_{i=1}^{r-1}ik_i\ge \frac{r}{2}(2h+|k_r|-2).
\end{equation}
Therefore $|k_r|-2\le 2g-2$ and $|k_r|\le 2g$.
Finally,
$\sum_{i=1}^{r-1}k_i=|k_r|+2-2h\le 2g+2.$
\end{proof}

\begin{remark}
The bounds on $r$ and on $k_i$ are sharp. Indeed, consider a hyperelliptic curve $P$ covering $\mathbb{CP}^1$ with ramifications at the vertices of a regular $(2g+1)$-gon and at its center. The covering can be chosen such that the automorphism of $P$ induced by the rotation of this polygon acts nontrivially in the fibers and hence has order $r=2(2g+1)=4g+2$.

On the other hand, consider a hyperelliptic curve $C$ with involution $\tau$. We have $$\chi(C_1(\tau))=2g+2, \chi(C_2(\tau))=2-2g-(2g+2)=-4g,$$
hence a pair $(C,\tau)$ belongs to the moduli space 
$\mathcal{M}_{g}(2g+2,-2g).$ 
\end{remark}

\begin{theorem} 
\label{th:EqToOrb}
The following equation holds:
\begin{equation}
\label{eq:EqToOrb}
\sum_{n=0}^{\infty}t^n\chi^{S_{n}}(\CM_{g,n})=\sum_{\underline{k}}\chi^{orb}(\CM_{g}(k_1,\ldots,k_r))\cdot\prod_{j=1}^{r}(1+p_jt^{j})^{k_j}.
\end{equation}
\end{theorem}

\begin{proof}
Consider the forgetful map $\pi_{g,n}:\CM_{g,n}\to \CM_{g}$. Its fiber over a point representing a curve $C$ is isomorphic to $F(C,n)/\Aut(C)$, hence one can apply Theorem \ref{th:conf} to compute its equivariant Euler characteristic:
$$
\sum_{n=0}^{\infty}t^n\chi^{S_n}(\pi_{g,n}^{-1}(C))=\sum_{n=0}^{\infty}t^n\chi^{S_n}(F(C,n)/\Aut(C))=$$ $$
\frac{1}{|\Aut(C)|}\sum_{\tau\in \Aut(C)}\prod_{i}(1+p_it^{i})^{\frac{\chi(C_{i}(\tau))}{i}}=
\sum_{\underline{k}}\frac{1}{|\Aut(C)|}\sum_{\tau\in \Aut_{\kk}(C)}\prod_{i}(1+p_it^{i})^{k_i}.
$$
Therefore:
$$
\sum_{n=0}^{\infty}t^n\chi^{S_{n}}(\CM_{g,n})=\int_{\CM_{g}}\sum_{n=0}^{\infty}t^n\chi^{S_n}(\pi_{g,n}^{-1}(C))d\chi=
$$
$$
\sum_{\underline{k}}\prod_{i}(1+p_it^{i})^{k_i}\int_{\CM_{g}}\frac{|\Aut_{\kk}(C)|}{|\Aut(C)|}d\chi.
$$
On the other hand, 
$$\chi^{orb}(\CM_{g}(k_1,\ldots,k_r))=\int_{\CM_{g}}\frac{|\pi_{g,\kk}^{-1}(C)|}{|\Aut(C)|}d\chi=\int_{\CM_{g}}\frac{|\Aut_{\kk}(C)|}{|\Aut(C)|}d\chi \qedhere
$$

\end{proof}
Using the Proposition \ref{prop:bounds}, we conclude that the sum in the right hand side of \eqref{eq:EqToOrb} is finite.

\begin{corollary}
The generating function $\sum_{n=0}^{\infty}t^n\chi^{S_{n}}(\CM_{g,n})$ is a rational function in $t$. Furthermore,  for any $n$,
$$
\chi^{S_{n}}(\CM_{g,n})\in \mathbb{Z}[p_1,\ldots,p_{4g+2}].
$$
\end{corollary}

The orbifold Euler characteristic of $\CM_{g}(k_1,\ldots,k_r)$ can be computed using the combinatorial results of Harer and Zagier \cite{harzag}.
We will denote the greatest common divisor of integers $a$ and $b$ by $(a,b)$. Let $\varphi(n)$ and $\mu(n)$ denote the Euler function and the M\"obius function respectively.
Define $$c(k,l,d):=\mu\left(\frac{d}{(d,l)}\right)\frac{\varphi(k/l)}{\varphi(d/(d,l))},$$

\begin{definition}
Let $\lambda=(\lambda_1,\ldots,\lambda_s)$ be a partition. We define a number
$$
N(r;\lambda)=\left|\left\{(x_1,\ldots,x_s)\in (\mathbb{Z}/r\mathbb{Z})^s : x_1+\ldots+x_s\equiv 0\ (\mbox{\rm mod}\ r), (x_i,k)=\lambda_i\right\}\right|.
$$
\end{definition}

\begin{lemma}(\cite{harzag})
\label{L2}
The following equation holds:
$$
N(r;\lambda)=\frac{1}{r}\sum_{d|r}\varphi(d)\prod_{i=1}^{s}c(k,\lambda_i,d).
$$
\end{lemma}

\begin{theorem}(\cite{harzag})
The orbifold Euler characteristic of the moduli space $\CM_{h,s}$ of genus $h$ curves with $s$ marked points is given by the formula:
$$
\chi^{orb}(\CM_{h,s})=(-1)^{s}\frac{(2h-1)B_{2h}}{(2h)!}(2h+s-3)!
$$
where $B_{k}$ denote Bernoulli numbers.
\end{theorem} 
 
\begin{theorem}
\label{fin}
The generating function for the $S_n$-equivariant Euler characteristics of $\CM_{g,n}$ has the form
$$
\sum_{n=0}^{\infty}t^{n}\chi^{S_n}(\CM_{g,n})=\sum_{\underline{k}}c_{k_1,\ldots,k_{r}}
\prod_{j=1}^{r}(1+p_jt^j)^{k_j},
$$
where $p_j$ are power sums  and the coefficients
$c_{k_1,\ldots,k_{r}}$ are defined  by the equation:
\begin{equation}
\label{defC}
c_{k_1,\ldots,k_r}=\chi^{orb}(\CM_{h,s})r^{2h}\prod_{p|\gamma}(1-p^{-2h})\cdot \frac{N(r;\lambda)}{ r\prod_{i=1}^{r-1}k_i!}.
\end{equation}
Here $h=\frac{1}{2}(1-\sum_{j=1}^{r}k_j), s=\sum_{j=1}^{r-1}k_j,$ $\gamma=\GCD(i:k_i>0)$ , $\lambda=\left(1^{k_1}2^{k_2}\ldots(r-1)^{k_{r-1}}\right)$
\end{theorem}

\begin{proof}
By Theorem \ref{th:EqToOrb} one has
$c_{k_1,\ldots,k_r}=\chi^{orb}(\CM_{g}(k_1,\ldots,k_r)).$
Consider the moduli space $\CM_{g}(k_1,\ldots,k_r)$ of pairs $(C,\tau)$.  As in Proposition \ref{prop:bounds}, to such a pair one can associate a genus $h$ curve $C_1=C/\tau$. The projection from $C$ to $C_1$ is ramified in $s$ points subdivided  into groups of size $k_1,\ldots,k_{n-1}$. The orbifold Euler characteristic of the moduli space of genus $h$ curves with such markings equals $\chi^{orb}(\CM_{h,s})/\prod_{i=1}^{r-1}k_i!$.

The number of pairs $(C,\tau)$ associated to a curve $C_1$ with fixed marked points was computed in \cite[pages 478--479]{harzag} and equals
$$
\frac{1}{r}r^{2h}\prod_{p|\gamma}(1-p^{-2h})\cdot N(r;\lambda).
$$
This completes the proof.
\end{proof}

The non-equivariant Euler characteristic of $\mathcal{M}_{g,n}$ has been computed in \cite[Theorem 4.3]{binihar}. It can be compared with Theorem \ref{fin} since 
$$\chi(\CM_{g,n})=n!\cdot \chi^{S_n}(\CM_{g,n})[p_1=1,\ p_k=0\ \mbox{\rm for}\ k>1].$$

\begin{example}
The generating function for the $S_n$-equivariant Euler characteristics of the moduli spaces of genus 2 curves with marked points has a form \cite{my}:

$$\sum_{n=0}^{\infty}t^n\chi^{S_n}(\CM_{2,n})=-\frac{1}{240}(1+p_1t)^{-2}-\frac{1}{240}(1+p_1t)^6(1+p_2t^2)^{-4}+$$
$$+\frac{2}{5}(1+p_1t)^3(1+p_5t^5)^{-1}+\frac{2}{5}(1+p_1t)(1+p_2t^2)(1+p_5t^5)(1+p_{10}t^{10})^{-1}+$$
$$+\frac{1}{6}(1+p_1t)^2(1+p_2t^2)(1+p_6t^6)^{-1}-\frac{1}{12}(1+p_1t)^4(1+p_3t^3)^{-2}-$$
$$-\frac{1}{12}(1+p_2t^2)^2(1+p_3t^3)^2(1+p_6t^6)^{-2}+\frac{1}{12}(1+p_1t)^2(1+p_2t^2)^{-2}+$$
$$+\frac{1}{4}(1+p_1t)^2(1+p_4t^4)(1+p_8t^8)^{-1}-\frac{1}{8}(1+p_1t)^2(1+p_2t^2)^2(1+p_4t^4)^{-2}.$$

These coefficients can be matched with the ones defined in Theorem \ref{fin}.
\end{example}

\end{document}